\documentclass{gtart_h}

\def\ifplaintex{\expandafter\ifx\csname documentclass\endcsname\relax}

\def\gtp{{\mathsurround=0pt\it $\cal G\mskip-2mu$eometry \&\ 
$\cal T\!\!$opology $\cal P\!$ublications}}  

\def\recd{{\small Received:\qua\receiveddate\ifx\reviseddate\relax
\else\qquad Revised:\qua\reviseddate\fi\par}} 

\def\lognumber#1{\def\thelognumber{#1}}
\def\volumenumber#1{\def\thevolumenumber{#1}}
\def\volumeyear#1{\def\thevolumeyear{#1}}
\def\papernumber#1{\def\thepapernumber{#1}}
\def\pagenumbers#1#2{\def\startpage{#1}\def\finishpage{#2}}
\def\published#1{\def\publishdate{#1}}

\def\received#1{\def\receiveddate{#1}}
\def\revised#1{\def\reviseddate{#1}}
\def\accepted#1{\def\accepteddate{#1}}

\long\def\asciiabstract#1{\long\def\theasciiabstract{#1}}

\let\\\par\let\thelognumber\relax\let\thevolumenumber\relax
\let\thepapernumber\relax\let\thevolumeyear\relax\let\startpage\relax
\let\finishpage\relax\let\publishdate\relax\let\receiveddate\relax
\let\reviseddate\relax\let\accepteddate\relax\let\theasciititle\relax
\let\theasciiauthors\relax
\let\theasciiabstract\relax

\let\theasciiemail\relax

\ifplaintex
\font\logobig=cmssbx10 scaled 3836
\font\logomed=cmssbx10 scaled 2557
\else
\font\logobig=cmssbx10 scaled 4200
\font\logomed=cmssbx10 scaled 2800
\fi

\long\def\makeagttitle{   
\count0=\startpage
\agt\hfill      
\hbox to 45truept{\vbox to 0pt{\vglue -13truept{\logomed A\kern -.37em{\logobig 
T}\kern -.38em G}\vss}\hss}
\break
{\small Volume \thevolumenumber\ (\thevolumeyear)
\startpage--\finishpage\nl
Published: \publishdate}

\vglue .25truein

{\parskip=0pt\leftskip 0pt plus
1fil\def\\{\par\smallskip}{\Large\bf\thetitle}\par\medskip} \vglue
0.05truein

{\parskip=0pt\leftskip 0pt plus 1fil\def\\{\par}{\sc\theauthors}
\par\medskip}%
 
\vglue 0.03truein 

{\small\leftskip 25truept\rightskip 25truept{\bf Abstract}\stdspace\theabstract

{\bf AMS Classification}\stdspace\theprimaryclass
\ifx\thesecondaryclass\relax\else; \thesecondaryclass\fi\par
{\bf Keywords}\stdspace \thekeywords\par}\vglue 7truept

}   

\ifplaintex
\font\phead=cmsl9 scaled 950
\font\pnum=cmbx10 scaled 913
\font\pfoot=cmsl9 scaled 950
\headline{\vbox to 0pt{\vskip -4.5mm\line{\small\phead\ifnum
\count0=\startpage ISSN 1472-2739 (on-line) 1472-2747 (printed)
\hfill {\pnum\folio}\else\ifodd\count0\def\\{ }%
\ifx\theshorttitle\relax\thetitle\else\theshorttitle\fi\hfill{\pnum\folio}
\else\def\\{ and }{\pnum\folio}\hfill\ifx\theshortauthors\relax\theauthors
\else\theshortauthors\fi\fi\fi}\vss}}
\footline{\vbox to 0pt{\vglue 0mm\line{\small\pfoot\ifnum\count0=\startpage
\copyright\ \gtp\hfill\else
\agt, Volume \thevolumenumber\ (\thevolumeyear)\hfill\fi}\vss}}
\else
\font=cmsl9 scaled 1050
\font\lnum=cmbx10 
\font=cmsl9 scaled 1050
\makeatletter
\def\@oddhead{{\small\lhead\ifnum\count0=\startpage ISSN 1472-2739 
(on-line) 1472-2747 (printed)\hfill {\lnum\number\count0}\else\ifodd\count0
\def\\{ }\ifx\theshorttitle\relax \thetitle \else\theshorttitle\fi\hfill
{\lnum\number\count0}\else\def\\{ and }{\lnum\number\count0}
\hfill\ifx\theshortauthors\relax 
\theauthors\else\theshortauthors\fi\fi\fi}}\def\@evenhead{\@oddhead}
\def\@oddfoot{\small\lfoot\ifnum\count0=\startpage\copyright\ \gtp\hfill\else
\agt, Volume \thevolumenumber\ (\thevolumeyear)\hfill\fi}
\def\@evenfoot{\@oddfoot}
\makeatother
\fi
\let\maketitlepage\makeagttitle
\let\makeshorttitle\maketitlepage
\let\maketitle\maketitlepage

\newwrite\gtoutfile
\long\gdef\makeheadfile{  
{\def\\{, }\def\s{ }
\immediate\openout\gtoutfile head.xxx
\immediate\write\gtoutfile{Proxy-for: \ifx\theasciiauthors\relax
\theauthors\else\theasciiauthors\fi\s<\ifx\theasciiemail\relax\theemail\else\theasciiemail\fi>}
\immediate\write\gtoutfile{\noexpand\\}
\immediate\write\gtoutfile{Authors: \ifx\theasciiauthors\relax
\theauthors\else\theasciiauthors\fi}
{\def\\{ }\immediate\write\gtoutfile{Title: \ifx\theasciititle\relax
\thetitle\else\theasciititle\fi}}
\immediate\write\gtoutfile{Subj-class: GT or SG, GR etc}
\immediate\write\gtoutfile{MSC-class: \theprimaryclass\ifx\thesecondaryclass\relax\else, \thesecondaryclass\fi}
\immediate\write\gtoutfile{Journal-ref: Algebr. Geom. Topol. \thevolumenumber\s
(\thevolumeyear) \startpage-\finishpage}
\immediate\write\gtoutfile{Comments: Published by Algebraic and
Geometric Topology at}
\immediate\write\gtoutfile{\s\s\s  http://www.maths.warwick.ac.uk/agt/AGTVol\thevolumenumber/agt-\thevolumenumber-\thepapernumber.abs.html}
\immediate\write\gtoutfile{\noexpand\\}
\immediate\write\gtoutfile{}
\ifx\theasciiabstract\relax
\immediate\write\gtoutfile{\theabstract}\else
\immediate\write\gtoutfile{\theasciiabstract}\fi
\immediate\write\gtoutfile{}
\immediate\write\gtoutfile{\noexpand\\}
\immediate\write\gtoutfile{}
\immediate\closeout\gtoutfile}}  

\def\maketitlepage{\makeagttitle\makeheadfile}
\let\makeshorttitle\maketitlepage
\let\maketitle\maketitlepage

\lognumber{62}
\volumenumber{5}
\volumeyear{2005}
\papernumber{62}
\pagenumbers{1555}{1572}
\received{15 December 2003} 
\revised{25 March 2005} 
\accepted{10 November 2005}
\published{23 November 2005}

\usepackage{amssymb,amsmath,amscd} 
\newtheorem{thm}{Theorem}[section]
\newtheorem{lem}[thm]{Lemma}		
\newtheorem{prop}[thm]{Proposition}
\theoremstyle{definition}
\newtheorem{defn}[thm]{Definition}
\newtheorem*{rem}{Remark}	
\newtheorem{example}[thm]{Example}
\newcommand{\smashprod}{\wedge}
\newcommand{\union}{\bigcup}
\newcommand{\RR}{{\mathbb R}}
\newcommand{\ls}[2]{\Omega^{#1}\Sigma^{#2}}
\newcommand{\ep}{\varepsilon}
\newcommand{\I}{{\mathcal I}}
\newcommand{\config}[2]{C(\RR^{#1}, #2)}
\newcommand{\intvl}[2]{I_{#1}(#2)}
\begin{document}

\title{The space of intervals in a Euclidean space}
\authors{Shingo Okuyama}				  
\address{Takuma National College of Technology\\Kagawa 
769-1192, JAPAN}
\email{okuyama@dc.takuma-ct.ac.jp}					   

\begin{abstract}   
For a path-connected space $X$, a well-known theorem of Segal, May and
Milgram asserts that the configuration space of finite points in
$\RR^n$ with labels in $X$ is weakly homotopy equivalent to
$\ls{n}{n}X$. In this paper, we introduce a space $\intvl{n}{X}$ of
intervals suitably topologized in $\RR^n$ with labels in a space $X$
and show that it is weakly homotopy equivalent to $\ls{n}{n}X$ without
the assumption on path-connectivity.
\end{abstract}

\asciiabstract{For a path-connected space X, a well-known theorem of
Segal, May and Milgram asserts that the configuration space of finite
points in R^n with labels in X is weakly homotopy equivalent to the
n-th loop-suspension of X.  In this paper, we introduce a space I_n(X)
of intervals suitably topologized in R^n with labels in a space X and
show that it is weakly homotopy equivalent to n-th loop-suspension of
X without the assumption on path-connectivity.}


\maketitle  

\section{Introduction}
G.Segal \cite{segal} introduced the configuration space
$\config{n}{X}$ of finite number of points in $\RR^n$ with labels in a
space $X$ and showed that $\config{n}{X}$ is weakly homotopy
equivalent to $\ls{n}{n}X$ if $X$ is path-connected.  When $X$ is not
path-connected, it follows from Segal's result that $\ls{n}{n}X$ is a
group-completion of $\config{n}{X}$, i.e.\ that
$H_*(\config{n}{X};k)[\pi^{-1}]$ is isomorphic to $H_*(\ls{n}{n}X;k)$
for any field $k$, where $[\pi^{-1}]$ denotes the localization of the
Pontrjagin ring $H_*(\config{n}{X};k)$ with respect to a sub-monoid
$\pi=\pi_0(\config{n}{X})$. (This was also shown independently by
F.Cohen \cite{cohen}.)  On the other hand, in \cite{mcduff}, D.McDuff
considered the space $C^{\pm}(M)$ of positive and negative particles
in a manifold $M$ and showed that it is weakly equivalent to some
space of vector fields on $M$. The topology of $C^{\pm}(M)$ is given
so that two particles cannot collide if they have the same parity, but
they can collide and annihilate if they are oppositely charged.  When
$M=\RR^n$, we can think of $C^{\pm}(\RR^n)$ as a $H$-space obtained by
adjoining homotopy inverses to $\config{n}{S^0}$.  Since adjoining
homotopy inverses to a $H$-space is a sort of group completion, one
might hope that $C^{\pm}(\RR^n)$ is weakly equivalent to
$\Omega^n\Sigma^n S^0=\Omega^nS^n$, but in fact,
$C^{\pm}(\RR^n)\simeq_w \Omega^n(S^n\times S^n/\Delta)$, where
$\Delta$ is the diagonal subspace \cite{mcduff}.  The aim of this
paper is to construct a configuration space model which is a
group-completion of $\config{n}{X},$ thus is weakly homotopy
equivalent to $\ls{n}{n}X$ for any $X$.

Caruso and Waner \cite{caruso-waner} constructed such a group-completion
model based on the space of little cubes \cite{may}.
They constructed the space of ``signed cubes merged along
the first coordinate" and showed that it approximates $\ls{n}{n}X$ without
the assumption on path-connectivity of $X$.
In this paper, we introduce a space $\intvl{n}{X}$ of intervals 
suitably topologized in $\RR^n$ and show that it gives another model for
the group-completion.
Our construction is, in some sense, a direct
generalization of $\config{n}{X}$ and simpler than the Caruso-Waner model.
More precisely, $\intvl{n}{X}$ is the space of intervals
ordered along parallel axes in $\RR^n$ with labels in $X$. In
this space the topology is such that ``cutting and pasting" and ``birth and death"
of intervals are allowed; i.e.\
cutting and pasting means that an interval with an open end and another
interval with a closed end can be attached at those ends to constitute 
one interval if they have the same label in $X$; birth and death means
that any half-open interval can vanish when its length tends to zero.

Now we can state our main theorem as follows:
\begin{thm} \label{thm:main}
	There is a weak homotopy equivalence $I_n(X)\simeq_w \ls{n}{n}X$.
\end{thm}

As contrasted with particles, intervals have two obvious features:
firstly, they are stretched in a direction, thus have a length; secondly,
any interval can be supposed to have a charge $(p,q)$ where
$p$ (resp.\ $q$) is $+1$ or $-1$  depending on
whether the interval contains or not contains the left (resp.\ right) endpoint.
The former feature results in the gradual interaction of our objects. 
For example, a possible process of annihilation of a closed and an open interval
is that: firstly, they are attached into one half-open interval and
then its length decreases and finally it vanishes --- i.e.\ they gradually
annihilate.
The latter feature of the interval plays an important role when we construct
an analogue of the electric field map \cite{segal} from a``thickened version'' of $I_n(X)$
into $\Omega\config{n-1}{X}$. The main step of our proof of Theorem \ref{thm:main} is 
first deforming $I_n(X)$ into this thickened but equivalent version, then constructing
the map in question and showing it is a weak equivalence using quasifibration techniques.
We then conclude using Segal's classical result as applied to $\Sigma X$ which is
now path-connected.

This paper can be considered as a first step of a larger project
proposed by K.Shimakawa. His idea is to use manifolds in $G$-vector
spaces to approximate $\Omega^{V^{\infty}}\Sigma^{V^{\infty}}X$
equivariantly, where $G$ is a compact Lie group and $V^\infty$ is
an orthogonal $G$-vector space which contains all the irreducible
$G$-representations infinitely many times as direct summands.


In \S 2, we settle the notation for the configuration space with
labels in a partial abelian monoid and observe some of its properties.
The definition of $I_n(X)$ is given in \S 3. In the same section,
we construct modifications of $I_n(X)$ which is needed to prove Theorem 
\ref{thm:main}. In \S 4, we construct a map
$\alpha\co \widetilde{I}_n(X)\to \Omega \config{n-1}{\Sigma X}$ and state
Proposition \ref{thm:quasifibration}, which is the key to prove Theorem 
\ref{thm:main}.
A proof of Theorem \ref{thm:main} is given in the same section.
We give a proof of Proposition \ref{thm:quasifibration} in \S 5.
Throughout the paper, $X$ is assumed to be a space with 
non-degenerate base point $*$.

The author is indebted to K.Shimakawa for many helpful suggestions.
He is also grateful to the referee for many valuable comments.
\section{Configuration space and partial abelian monoid} \label{sec:functor}
The notion of the configuration space with summable labels appear
in several papers \cite{salvatore},\cite{kallel},\cite{shimakawa}.
In this section, we introduce one form of such notion adapted
to the purpose of this paper. Remark that our definition of
partial abelian monoid is a special case of that given in \cite{shimakawa},
of which we shall call `two-generated' partial abelian monoid.

\begin{defn}
A partial abelian monoid (PAM for short) is a space $M$ equipped with
a subspace $M_2\subset M\times M$ and a map $\mu\co M_2\to M$ such that
\begin{enumerate}
	\item $M\vee M\subset M_2$, $\mu(a,*_M)=\mu(*_M,a)=a$,	
	\item $(a,b)\in M_2\Leftrightarrow (b,a)\in M_2, \mu(a,b)=\mu(b,a)$, and
	\item $(\mu(a,b),c)\in M_2\Leftrightarrow (a,\mu(b,c))\in M_2,
	\mu(\mu(a,b),c)=\mu(a,\mu(b,c))$.
\end{enumerate}
We write $\mu(a,b)=a+b$. An element in $M_2$ is called a summable pair.
Let $M_k$ denote the subspace of $M^k$ which consists of those
$k$-tuples $(a_1,\dots,a_k)$ such that $a_1+\dots +a_k$ is defined.
A map between PAMs are called a PAM homomorphism if it 
sends summable pairs to summable pairs and preserves the sum.
\end{defn}

\begin{defn}
	Let $Z^{(k)}_n(M)$ denote the subspace of
$(\RR^n\times M)^k$ given by:
\[
	Z^{(k)}_n(M)=\left\{((v_1,a_1),\dots,(v_k,a_k))\left|%
	\begin{array}{l}\mbox{for any~} i_1,\dots, i_r
	\mbox{~such that~}\\\hspace{2cm}v_{i_1}=\dots=v_{i_r},\\
	(a_{i_1},\dots, a_{i_r})\in M_r\\
	\end{array}\right.\right\}
\]
Then we define a space $\config{n}{M}$ as $\config{n}{M}=(\coprod_{k\geq 0}Z^{(k)}_n(M))/\sim$,
where $\sim$ denotes the least equivalence relation which satisfies 
(R1)$\sim$(R3) below.
\begin{enumerate}
	\item[(R1)] If $a_i=*_M$ then
	\[
		((v_1,a_1),\dots,(v_k,a_k))%
		\sim ((v_1,a_1),\dots,\widehat{(v_i,a_i)},\dots,(v_k,a_k)).
	\]
	\item[(R2)] For any permutation $\sigma\in \Sigma_k$,
	\[
		((v_1,a_1),\dots, (v_k,a_k))%
		\sim ((v_{\sigma^{-1}(1)},a_{\sigma^{-1}(1)}),\dots,%
		(v_{\sigma^{-1}(k)},a_{\sigma^{-1}(k)})).%
	\]
	\item[(R3)] If $v_1=v_2=v$,
	\[
		((v_1,a_1),\dots,(v_k,a_k))%
		\sim ((v,a_1+a_2),(v_3,a_3),\dots,(v_k,a_k)).%
	\]
\end{enumerate}

We regard $\config{n}{M}$ as a PAM as follows.
Let $\config{n}{M}_2\subset \config{n}{M}^2$ denote the subspace which consists of pairs
$(\xi, \eta)$ 
which have representatives
$((v_1,a_1),\dots,(v_k,a_k))$ and 
$((v_1,b_1),\dots,(v_k,b_k))$ 
such that $v_i$'s are distinct and $(a_i,b_i)\in M_2$
for every $i$.
(Note that some of $a_i$'s or $b_j$'s may be zero.)
 Then $\mu\co\config{n}{M}_2\to \config{n}{M}$ is defined by setting
$\mu(\xi,\eta)=[(v_1,a_1+b_1),\dots,(v_k,a_k+b_k)]$.
Note that $\config{n}{-}$ is a self-functor on
the category of partial abelian monoids.
\end{defn}

\begin{example}
	Any space $X$ is regarded as a PAM by
	setting $X_2=X\vee X$. Then $\config{n}{X}$ is nothing but the 
	configuration
	space of finite points in $\RR^n$ labelled by $X$.
\end{example}

\begin{example} \label{ex:mcduff}
	Let $M=\{-1,0,1\}$ with $1+(-1)=0$ as the only non-trivial partial	
	sum. Then $\config{n}{M}$ is homeomorphic to $C^\pm(\RR^n)$, the configuration space of
	positive and negative particles in $\RR^n$, given in \cite{mcduff}.
	Furthermore, we give $X\wedge M$ a PAM structure by setting
	$(X\wedge M)_2=\{(x,a;x,b)~|~(a,b)\in M_2\}$ and $(x,a)+(x,b)=(x,a+b)$. Then
	$\config{n-1}{X\wedge M}$ is the labelled version of $C^{\pm}(\RR^n)$ \cite{caruso-waner}.
\end{example}

Lemma \ref{thm:htp_functor} below states that the functor $\config{n}{-}$
preserves homotopies. Homotopy in the category of PAMs is
defined as follows:
Let $M$ be a PAM. We regard $M\times I$ as a PAM by setting
$$
(M\times I)_2=\{(m,t;n,t) ~|~ (m,n)\in M_2\}
$$
$$\mu_{M\times I}(m,t;n,t)=(\mu_M(m,n),t).
\leqno{\hbox{and}}
$$
Homomorphisms $f,g\co M\to N$ between PAMs are
called homotopic via PAM homomorphisms if
there exists a homomorphism $H\co M\times I\to N$ such that
$H_0=f, H_1=g$.

\begin{lem} \label{thm:htp_functor}
	If $f,g\co M\to N$ are homotopic via PAM homomorphisms then
	$\config{n}{f}$ and $\config{n}{g}$ are homotopic via PAM homomorphisms.
\end{lem}
\begin{proof}
	Let $H\co M\times I\to N$ be a homotopy between $f$ and $g$.
	Observe that we have a homomorphism 
	$\config{n}{M}\times I\to \config{n}{M\times I}$
	by setting	
	\[
		([(v_1,a_1),\dots,(v_k,a_k)],t)\mapsto%
		[(v_1,(a_1,t)),\dots,(v_k,(a_k,t))].
	\]
	Then this map followed by $\config{n}{H}$ is a homotopy between
	$\config{n}{f}$ and $\config{n}{g}$ via PAM homomorphisms.
\end{proof}

Recall that the Moore loop space on a space $X$ is defined by
$$
	\Omega (X)=\cup_{s\geq 0}\Omega_s(X)\times \{s\},
$$
$$
	\Omega_s(X)=\{l\co [0,s]\to X ~|~ l(0)=l(s)=*\}
\leqno{\hbox{where}}
$$
is the space of loops of length $s$.
Recall also that $\Omega(X)$ is
topologized as the subspace of $\mbox{Map}([0,\infty), X)\times [0, \infty)$.
Let $M$ be a PAM. We give $\Omega_s(M)$ a PAM structure
by setting
$$
	(\Omega_s(M))_2=%
	\{(l_1,l_2) ~|~ (l_1(t),l_2(t))\in M_2 \mbox{~for all~}t\in [0,s]\}
$$
$$(l_1 + l_2)(t) = l_1(t) + l_2(t) \in \Omega_s(M),\ t \in [0,s] .
\leqno{\hbox{and}}
$$
It is clear that $\Omega_s$ is a self-functor on the category of PAMs.
We have a map $\config{n}{\Omega_s (M)}\to\Omega_s\config{n}{M}$ defined by
\[
	[(v_1,l_1),\dots, (v_k,l_k)]\mapsto %
	(t\mapsto [(v_1,l_1(t)),\dots,(v_k,l_k(t))]),
\]
which will be used in the construction of the map 
$\alpha\co \widetilde{I}_n(X)\to \Omega \config{n-1}{\Sigma X}$
in \S \ref{sec:map}. 

\section{The space of intervals in $\RR^n$}
\label{sec:definition}

Let $\I$ denote the subspace of $\RR^2\times \{\pm 1\}^2$ consisting of
all quadruples $(u,v,p,q)$ such that $u<v$ if $p=q$ and $u\leq v$ if $p\neq q$.
When $u<v,$ $J=(u,v,p,q)\in \I$ can be identified with an interval in $\RR$ whose
endpoints are $u$ and $v$.
It contains (resp.\ not contains) the left endpoint if $p=1$ (resp.\ $p=-1$).
Similarly, it contains (resp.\ not contains) the right endpoint if
$q=1$ (resp.\ $q=-1$).
Thus, for example, $J=(u,v, -1, 1)$ with
$u<v$ is identified with the half-open interval $(u,v]$.
For any $J=(u,v,p,q)\in \I,$ we put $l(J)=u, r(J)=v, p_L(J)=p$ and
$p_R(J)=q$. 
If $J,K\in \I$ and $r(J)\leq l(K)$ then we write $J\leq K$. If, moreover, $r(J)<l(K)$
then we write $J<K$. 

Let $U$ be a connected subset of $\RR$. 
In applications, $U$ is one of $(-s,s)$ or $(0,s)$ for $s>0$ or $s=\infty$.
Let $\I(U)$ denote the subspace of $\I$ consisting of those $J\in \I$
such that $l(J),r(J)\in U$. Remark that all the intervals in
$\I=\I(\RR)$ are bounded.

\begin{defn}\label{defn:I_(k)X_U}
We define $I_{(k)}(X)_U$ as the subspace of
$(\I(U)\times X)^k$ consisting of those $k$-tuples
$((J_1,x_1),\dots, (J_k, x_k))$ such that
\begin{enumerate}
	\item[(1)] $J_1\leq J_2\leq \dots \leq J_k$,
	\item[(2)] $x_{i-1}\neq x_i$ implies $J_{i-1}< J_i$, and
	\item[(3)] $p_R(J_{i-1}) = p_L(J_i)$ implies  $J_{i-1}<J_i$.
\end{enumerate}
In other words, $J_1, \dots, J_k$ are disjoint ordered intervals
of $U$ with labels\break $x_1, \dots, x_k$ and
$J_{i-1}$ and $J_i$ can have common endpoints
only if the labels in $X$ coincide and the given endpoints are of
the opposite sign.
\end{defn}

\begin{defn}\label{def:I1}
We define the space of intervals in $U$ to be
\[
	I_1(X)_U=\coprod_{k\geq 0}I_{(k)}(X)_U/\sim,
\]
where $\sim$ denotes the equivalence relation generated by the relation
shown below.
Suppose 
$$
	\iota=((J_1,x_1),\dots,(J_k,x_k))\in I_{(k)}(X)_U
$$
$$
	\iota'=((K_1,y_1),\dots,(K_{k-1},y_{k-1}))\in I_{(k-1)}(X)_U.
\leqno{\hbox{and}}
$$
Then $\iota'\sim\iota$ if one of the following holds:
\begin{enumerate}
	\item (cutting and pasting)
		\[
		K_i=\left\{\begin{array}{ll}
			J_i & \mbox{~if~} i<j\\
			J_j\cup J_{j+1} & \mbox{~if~} i=j\\
			J_{i+1} & \mbox{~if~} i>j ~,
		\end{array}\right.
		y_i=\left\{\begin{array}{ll}
			x_i & \mbox{~if~} i<j\\
			x_j=x_{j+1} & \mbox{~if~} i=j\\
			x_{i+1} & \mbox{~if~} i>j.
		\end{array}\right.
		\]
\smallskip
	\item (birth and death)
		\[
		K_i=\left\{\begin{array}{ll}
			J_i & \mbox{~if~} i<j\\
			J_{i+1} & \mbox{~if~} i\geq j ~,
		\end{array}\right.
		y_i=\left\{\begin{array}{ll}
			x_i & \mbox{~if~} i<j\\
			x_{i+1} & \mbox{~if~} i\geq j,
		\end{array}\right.\mbox{~and~}
		\]
		$x_j=*$ or $J_j=(u,u,p,-p)$ for some $u$ and $p$.
\end{enumerate}

For any $\iota\in I_1(X)_U,$
we have a representative $((J_1,x_1),\dots,(J_k,x_k))$ such that
$x_i\neq *$ for every $i$ and $J_1<J_2<\dots<J_k$, which is called
the reduced representative.

We regard $I_1(X)_U$ as a partial monoid by considering
a pair $(\xi,\eta)\in I_1(X)_U^2$ is summable if
$\xi$ and $\eta$ have representatives such that the union of
those
satisfies the conditions (1)$\sim$(3) in Definition 
\ref{defn:I_(k)X_U} after an appropriate change of the
order of labelled intervals.
The sum is given by the union of such representatives.
It is clear that the only element in $I_{(0)}(X)_U$, denoted $\emptyset$,
is the unit for the partial sum.
\end{defn}
\begin{rem}
	Recall from Example \ref{ex:mcduff} that when $M=\{-1,0,1\},$ 
	$\config{n}{X\wedge M}$ is the configuration space of positive and
	negative particles in $\RR^n$ labelled by $X$. Remark that $I_1(X)$ 
	can be embedded into $\config{}{X\wedge M}$ as a topological space in the following
	manner.
	Let $C'(\RR,X\smashprod M)$ denote
	the sub-PAM of $\config{}{X\smashprod M}$ consisting of elements
	in $\config{}{X\smashprod M}$ which have a representative
	\[
		((v_1,x_1\smashprod a_1),\dots,(v_{2k},x_{2k}\smashprod a_{2k}))%
	\in Z^{(2k)}_1(X\smashprod M)
	\]
	such that $v_1\leq \dots \leq v_{2k}$ and
	$x_{2i-1}=x_{2i}$ for all $i$.
	Then we have a homeomorphism $I_1(X)\to C'(\RR,X\smashprod M)$ defined by the
	correspondence
	\[
	[(J_1,x_1),\dots, (J_k,x_k)]\mapsto[(u_1,p_1,x_1),(v_1,q_1,x_1),\dots,(u_k,p_k,x_k),(v_k,q_k,x_k)],
	\]
	where $J_i=(u_i, v_i, p_i,q_i)$.
	We do not have such a relation between $I_n(X)$ and $\config{n}{X\wedge M}$
	for $n>1;$ the above homeomorphism does not help us about this 
	since it is not a PAM homomorphism.
\end{rem}

\begin{defn}
We define $I_n(X)_U=\config{n-1}{I_1(X)_U}$.
We denote $I_n(X)=I_n(X)_{\RR}$ and $I_n(X)_s=I_n(X)_{(0,s)}$. Then
$I_n(X)$ is homeomorphic to $I_n(X)_s$ for any $s>0$.
An element $[(v_1,\xi_1),\dots,(v_k,\xi_k)]\in I_n(X)=\config{n-1}{I_1(X)}$
can be thought of as intervals ordered along the lines parallel to the $x_1$-axis
through $v_i$ if we view $\RR^{n-1}$ as the $(x_2,\dots,x_n)$ hyperplane
in $\RR^n$.
Thus we call $I_n(X)$ the space of intervals in $\RR^n$.
\end{defn}

To relate $I_n(X)$ with $\ls{n}{n}X,$ we construct an analogue of the
electric field map in \cite{segal}.
However, there is no direct analogue of the electric field map on $I_n(X)$ and
we need a thickening of $I_n(X)$. We also need to modify $I_n(X)$
to get a space corresponding to the Moore loop space for the quasifibration
argument given in Theorem \ref{thm:main} and the related lemmas.

\begin{defn} \label{defn:separated}
We say that $\iota\in I_{(k)}(X)_U$ is $\ep$-separated if
\begin{enumerate}
	\item[(1)] $\iota\in I_{(k)}(X)_{Int(U-U^c_{\ep/2})}$,	
	where $U^c_{\ep/2}$ denotes the $\ep/2$-neighborhood 
	of the complement of $U$,
	\item[(2)] any two ends (of the same or distinct intervals) with
	the same parity are distant by at least $\ep$, and
	\item[(3)] any two intervals with the distinct labels in $X$ are distant by at 
	least $\ep$.
\end{enumerate}
We then say that $\xi\in I_1(X)_U$ is $\ep$-separated if
it is represented by some $\ep$-separated $\iota\in I_{(k)}(X)_U$.
\end{defn}

Let $I^\ep_1(X)_U$ denote the subspace of $I_1(X)_U$ consisting of 
$\ep$-separated elements.
Then $I^\ep_1(X)_U$ is given a PAM structure by regarding
$(\xi,\eta)\in (I^\ep_1(X)_U)^2$ as a summable pair if it is so as a pair
of elements in $I_1(X)_U$ and the sum $\xi+\eta$ taken there is in 
$I^\ep_1(X)_U$.
We define
\[
	I^\ep_n(X)_U=\config{n-1}{I^\ep_1(X)_U},
\]
then define
\[
	\widetilde{I}_n(X)=\{(\xi, \ep, s) ~|~ 0<\ep\leq 1,~s\geq 0,~%
	\xi\in I^\ep_n(X)_s \}%
\]
with the topology considered as the subspace of
$I_n(X)_\infty\times (0,1]\times [0,\infty)$.
A PAM structure on $\widetilde{I}_n(X)$ is defined so that
$(\xi,\ep, s)$ and	$(\eta,\tau, t)$ are summable if and only if
$\ep=\tau, s=t$ and  $(\xi, \eta)$ is a summable pair.

\begin{lem} \label{thm:thickening}
	We have a weak homotopy equivalence $\widetilde{I}_n(X)\simeq_w I_n(X)$.
\end{lem}

\begin{proof}
	Since any homeomorphism $\RR\to (0,\infty)$ induces
	a homeomorphism\break $I_n(X)\to I_n(X)_\infty$,
	it suffices to show that
	$\widetilde{I}_n(X)\simeq_w  I_n(X)_\infty$.
	Note that we can embed $I^\ep_n(X)_s$ into $I_n(X)_\infty$ using
	the inclusion $(0,s)\subset (0,\infty)$. Under this identification,
	let $p\co \widetilde{I}_n(X)\to I_n(X)_\infty$ denote the map
	which assigns $\xi$ to each $(\xi,\ep,s)$.
	Then $p_*\co \pi_k(\widetilde{I}_n(X))\to \pi_k(I_n(X)_\infty)$
	is an isomorphism for all $k\geq 0$.
	Indeed, for any map $f\co S^k\to I_n(X)_\infty$, there exist $\ep$ and
	$s$ such that $\mbox{Im} f$ is contained in $I^\ep_n(X)_s$ since $S^k$
	is compact. This proves that $p_*$ is surjective.
	On the other hand, let $H\co S^k\times I\to I_n(X)_\infty$ be a homotopy
	between $H_0=p\circ f$ and $H_1=p\circ g$. Since $S^k\times I$ is
	compact, we can restrict the codomain of $H$ to $I^\ep_n(X)_s$ for
	some $\ep$ and $s$. This proves that $p_*$ is injective.
\end{proof}

Now we proceed to construct another modification of $I_n(X),$
which models, as we shall see later, the Moore path space $P\config{n-1}{\Sigma X}$.
Let $s$ be a positive number, or $s=\infty$. 
For any element $J=(u,v,p,q)\in \I$, we put
$-J=(-v,-u,-q,-p)\in \I$.
Then we have an involution on
$I_{(k)}(X)_{(-s,s)}$ by setting
\[
	(-1)\cdot ((J_1,x_1),\dots, (J_k,x_k))=((-J_k,x_k),\dots, (-J_1, x_1)),
\]
which induces an involution on $I_1(X)_{(-s,s)}$.
We denote by $E_1(X)_s$, the subspace of
$I_1(X)_{(-s,s)}$ invariant under the involution. Note that $E_1(X)_s$ has 
a PAM structure induced by that of $I_1(X)_{(-s,s)}$;
we define $E_n(X)_s=\config{n-1}{E_1(X)_s}.$
Since the involution on $I_1(X)_{(-s,s)}$ restricts to an involution on
$I^{\ep}_n(X)_{(-s,s)}$, we can define $E^\ep_1(X)_s$ and
$E^\ep_n(X)_s=\config{n-1}{E^\ep_1(X)_s}$ similarly.
Now we define,
\[
	\widetilde{E}_n(X)=\{(\xi, \ep, s)~|~0< \ep\leq 1, s\geq 0,%
	\xi \in E^\ep_n(X)_s\},
\]
with the topology considered as a subset of
$E_n(X)\times (0, 1]\times [0, \infty)$.
A PAM structure on $\widetilde{E}_n(X)$ is defined so that
$(\xi,\ep, s)$ and	$(\eta,\tau, t)$ are summable if and only if
$\ep=\tau, s=t$ and  $(\xi, \eta)$ is a summable pair.

Any element of $E_1(X)$ have a representative of the form
\[
	((-J_k,x_k),\dots,(-J_1,x_1), (J_1,x_1),\dots,(J_k,x_k))%
	\in I_{(2k)}(X)_{(-s,s)}
\]
for some $s>0$.

It is useful to denote this representative by $m((J_1,x_1),\dots,(J_k,x_k))$.

\begin{lem} \label{thm:contraction}
	We have a weak homotopy equivalence
	$\widetilde{E}_n(X)\simeq_w *$.
\end{lem}
\begin{proof}
	We can prove that $\widetilde{E}_n(X)\simeq_w E_n(X)$ in a
	similar way to the proof of Lemma \ref{thm:thickening}.	
	So, by Lemma \ref{thm:htp_functor}, it suffices to show that
	$E_1(X)$ is homotopy equivalent to $\{0\}$ via PAM homomorphisms.
	Since $E_1(X)$ is homeomorphic to $E_1(X)_s$ for any $s>0$, we prove
	that $E_1(X)_1\simeq \{0\}$ via PAM homomorphisms.
	Let $h_t\co (-1,1)\to (-1,1)~(0\leq t\leq 1)$ denote the homotopy defined by:
	\[
		h_t(u)=\left\{\begin{array}{cc}
			u-t & \mbox{~if~} u\geq t\\
			0  & \mbox{~if~} |u|<t\\
			u+t & \mbox{~if~} u\leq -t
		\end{array}\right.
	\]
	Then the contracting homotopy $H_t\co E_1(X)_1\to E_1(X)_1$ is defined as follows.
	For $\xi\in E_1(X)$, we take a representative	
	\[
		m((J_1,x_1),\dots,(J_k,x_k))\in I_{(2k)}(X)_{(-1,1)}.
	\]
	Then $H_t\co E_1(X)_1\to E_1(X)_1$ is defined by setting
	$H_t(\xi)$ as the class represented by
	\[
		m((h_t(J_r),x_r),\dots,(h_t(J_k), x_k))\in I_{(2(k-r+1))}(X)_{(-1,1)},
	\]
	where $r$ is the least integer among $i>0$ such that $r(J_i)>t$ and
	$h_t(J_i)$ denotes the element $(h_t(l(J_i)),h_t(r(J_i)), p_L(J_i), p_R(J_i))\in \I$.
	It is straightforward to show that $H_t$ is the desired contracting homotopy.
\end{proof}

Consider the map $I_{(k)}(X)_s\to I_{(k)}(X)_{(-s,s)}$ given by
\[
	((J_1,x_1),\dots,(J_k,x_k))\mapsto %
	((-J_k,x_k),\dots,(-J_1,x_1),(J_1,x_1),\dots,(J_k,x_k)).
\]
These maps for all $k$ induce a map $I_1(X)_s\to E_1(X)_s$, 
which restricts to
a map $I^{\ep}_1(X)_s\to E^{\ep}_1(X)_s$.
Thus we have an embedding
$i\co \widetilde{I}_n(X)\to \widetilde{E}_n(X)$.
We regard $\widetilde{I}_n(X)$ as the subspace of $\widetilde{E}_n(X)$ via this
embedding.
\section{The map $\alpha\co \widetilde{I}_n(X)\to \Omega \config{n-1}{\Sigma X}$
and the proof of Theorem \ref{thm:main}} \label{sec:map}
We first construct a map $\alpha\co \widetilde{I}_n(X)\to \Omega \config{n-1}{\Sigma X}$.
Let $\xi$ be an element of $I^\ep_1(X)_s$ and
$((J_1,x_1), \dots, (J_k, x_k))\in I_{(k)}(X)_s$ denote
the reduced representative of $\xi$.
Let $u_{2i-1}=l(J_i), u_{2i}=r(J_i),
p_{2i-1}=p_L(J_i)$ and $p_{2i}=p_R(J_i)$.
Then 
we define $N_i\subset [0,s]~(i=1,\dots, 2k)$ as
\[
	N_1=\left[u_1-\ep/2, \mbox{Min}\left(u_1+\ep/2, u_{2}-\ep/2\right) \right],
\]
\[
	N_i=\left[\mbox{Max}\left(u_i-\ep/2, u_{i-1}+\ep/2\right), 
	\mbox{Min}\left(u_i+\ep/2, u_{i+1}-\ep/2\right) \right],	
	\mbox{~for~} 1<i<2k,
\]
$$	N_{2k}=\left[\mbox{Max}\left(u_{2k}-\ep/2, u_{2k-1}+\ep/2\right), u_{2k}+\ep/2
	\right].\leqno{\hbox{and}}$$

\begin{lem}
There exists a map $f\co [0,s]\to S^1\wedge X$ such that
\begin{enumerate}
\item $f(t)=[p_i(\frac{t-u_i}{\ep}+\frac{(-1)^i}{2})]\wedge x_{[\frac{i+1}{2}]}, \mbox{~if~} t\in N_i,$
where $[\frac{i+1}{2}]$ denotes the largest integer which does not exceed $\frac{i+1}{2},$
\item $f$ is piecewise constant outside $\union_{i=1}^{2k}N_i,$ and
\item $f(0)=f(s)=*,$ the base point of $S^1\wedge X$.
\end{enumerate}
\end{lem}
\begin{proof}
First of all, $N_i$'s are non-empty. The only problem is to show that
$u_{i-1}+\ep/2 < u_{i+1}-\ep/2$. However, at least two of three points located at
$u_{i-1}, u_i$ and $u_{i+1}$ should have the same parity. Hence we have 
$u_{i+1}-u_{i-1}\geq \ep$ by (2) of Definition \ref{defn:separated}.
Since we took the reduced representative of $\xi,$
it also follows that $N_i$'s are intervals such that $N_1\leq \dots \leq N_{2k}$.

By the definition of $N_i$'s, $N_i\cap N_{i+1}\neq\emptyset$ only
if $r(N_i)=l(N_{i+1})$ i.e.\ $u_{i+1}-u_i=\ep$. So suppose $u_{i+1}-u_i=\ep$.
If $i$ is even, then we have $f(r(N_i))=f(l(N_{i+1}))=*\in S^1\wedge X$.
On the other hand, if $i$ is odd, then we have $f(r(N_i))=f(l(N_{i+1}))=0\wedge x_{\frac{i+1}{2}}$.
Thus $f$ is well-defined on $\union^{2k}_{k=1}N_i$.

Next, we show that $f:\union^{2k}_{k=1}N_i\to S^1\wedge X$ can be extended to
$[0,s]$ so that it is piecewise constant outside $\union_{i=1}^{2k}N_i$.
To show this, it suffices to show that $f(r(N_i))=f(l(N_{i+1}))$ even for
$u_{i+1}-u_i\neq \ep$. If $u_{i+1}-u_i>\ep,$ 
$f(r(N_i))=f(l(N_{i+1}))$ follows by the same argument as above.
Suppose $u_{i+1}-u_i< \ep$. By (2) of Definition \ref{defn:separated}, we have
$p_i= -p_{i+1}$. We also have $x_{[\frac{i+1}{2}]}=x_{[\frac{i+2}{2}]}$,
by (3) of Definition \ref{defn:separated} if $i$ is even.
Now a direct calculation shows that $f(r(N_i))=f(l(N_{i+1}))$.

Finally $f(0)=f(s)=*$ since, by (1) of Definition \ref{defn:separated}, we have
$l(N_1)\geq 0$ and $r(N_{2k})\leq s$.
\end{proof}

By setting $\alpha^\ep_s(\xi)=f$, we obtain a PAM homomorphism
\[
		\alpha^\ep_s\co I^\ep_1(X)_s\to \Omega_s(\Sigma X),
\]
where $\Omega_s$ is the ``loops of length $s$" functor defined in \S\ref{sec:functor}.
Then we define a map
$\alpha\co \widetilde{I}_1(X)\to \Omega\Sigma X$ by 
$(\xi,\ep,s)\mapsto \alpha^\ep_s (\xi)$, which is also a
PAM homomorphism.
Now we define a map $\alpha\co \widetilde{I}_n(X)\to \Omega \config{n-1}{\Sigma X}$
by the composite
\[
	\widetilde{I}_n(X)\to \config{n-1}{\widetilde{I}_1(X)}\stackrel{\config{n-1}{\alpha}}%
	{\longrightarrow} \config{n-1}{\Omega\Sigma X}\to \Omega \config{n-1}{\Sigma X}.
\]
where the first map is defined similarly to the map given in the proof
of Lemma \ref{thm:htp_functor} and the last map is the one given in the
end of \S \ref{sec:functor}.

Next, we construct a map $p\co \widetilde{E}_n(X)\to \config{n-1}{\Sigma X}$.
We define
$$\alpha^\ep_{(-s, s)}\co I^\ep_1(X)_{(-s, s)}\to \mbox{Map}((-s,s);\Sigma X)$$
similarly to the definition of $\alpha^\ep_s$.
Let $p^\ep_s$ denote the composite
\[
	E^\ep_1(X)_s\rightarrow I^\ep_1(X)_{(-s, s)}%
	\stackrel{\alpha^\ep_{ (-s, s)}}{\longrightarrow}%
	\mbox{Map}((-s, s);\Sigma X)\stackrel{e_0}{\longrightarrow}\Sigma X,
\]
where $e_0$ denotes the evaluation at $0\in (-s, s)$.
The map
$p^\ep_s$ induces a map
\[
	p\co\widetilde{E}_n(X)\to \config{n-1}{\Sigma X},
\]
which is surjective.

\begin{rem} \label{rem:scan at zero}
By definition, a configuration $p(\xi,\ep, s)\in \config{n-1}{\Sigma X}$ has
a particle located at $v\in \RR^{n-1}$ with non-trivial label in $X$
if and only if the configuration $\xi\in E_n^\ep(X)_s=\config{n-1}{E^\ep_1(X)_s}$
has a particle located at $v$ labelled by some
$\iota=[m((J_1,x_1),\dots,(J_k,x_k))]\in E^\ep_1(X)_s$ such that
$l(J_1)<\ep$. Thus the fiber of $p$ at $\emptyset\in \config{n-1}{\Sigma X}$ is 
$\widetilde{I}_n(X)$.
\end{rem}

\begin{prop} \label{thm:quasifibration}
	The sequence
	$\widetilde{I}_n(X)\stackrel{i}{\to}\widetilde{E}_n(X)\stackrel{p}{\to} \config{n-1}{\Sigma X}$
	is a quasifibration.
\end{prop}

Proof of Proposition \ref{thm:quasifibration} is given in the next 
section.

\begin{proof}[Proof of Theorem 1]
We define $\beta\co\widetilde{E}_n(X)\to P\config{n-1}{\Sigma X}$ by
\[
		\beta(\xi,\ep,s)=\alpha^\ep_{(-s,s)}(\xi)|_{[0,s]}.
\]
By the definition of $\alpha$ and $\beta$ the following diagram is 
commutative.
\begin{equation*}
	\begin{CD}
		\widetilde{I}_n(X) @>{i}>>
			\widetilde{E}_n(X) @>{p}>>
				\config{n-1}{\Sigma X}\\
		@V{\alpha}VV
			@V{\beta}VV
				@|\\
		\Omega \config{n-1}{\Sigma X} @>>>
			P\config{n-1}{\Sigma X} @>>>
				\config{n-1}{\Sigma X},
	\end{CD}
\end{equation*}
where the lower horizontal row is the path-loop fibration on $\config{n-1}{\Sigma X}$.
Since $\beta$ is a weak homotopy equivalence by Lemma 
\ref{thm:contraction}, so is $\alpha$.
Now the theorem follows from the Segal's theorem \cite{segal} since $\Sigma X$
is path-connected.
\end{proof}
\section{Proof of Proposition \ref{thm:quasifibration}} 
\label{sec:prf_of_key}

Before stating the proof of Proposition \ref{thm:quasifibration}, we need some
observations on the filtration.

\begin{defn}
Let $M$ be a PAM and $A$ a closed sub-PAM of $M$.
The filtration on $\config{n}{M}$ associated to $A$ is
\[
	F^A_j\config{n}{M}=\{[(v_1,a_1),\dots,(v_k,a_k)] ~|~ \#\{i|a_i\notin A\}\leq j\}~(j\geq 0).
\]
\end{defn}

\begin{example}\label{ex:filt1}
	We have ``the number of points" filtration on $\config{n-1}{\Sigma X}$
	defined by
	\[
	F_j\config{n-1}{\Sigma X}=\coprod_{k\leq j} Z^{(k)}_{n-1}(\Sigma X)/\sim.
	\]
	This coincides with $F_j^A\config{n-1}{\Sigma X}$ if we put $A=*$.
\end{example}

\begin{example}\label{ex:filt2}
	Let $A^\ep_s$ denote the closed sub-PAM of $E^\ep_1(X)_s$ given by
	\[
		A^\ep_s=\{[m((J_1,x_1),\dots,(J_k,x_k))]\in%
		E^\ep_1(X)_s ~|~
		l(J_1)\geq \ep/2\}.
	\]
	Then we have a filtration $F^{A^\ep_s}_j\config{n-1}{E^\ep_1(X)_s}$ on
	$E^\ep_n(X)_s{=}\config{n-1}{E^\ep_1(X)_s}.$
	Using this filtration, we get a filtration on $\widetilde{E}_n(X)$
	by setting
	\[
		F_j\widetilde{E}_n(X)=\left\{ (\xi,\ep,s)~|~\xi\in
		F^{A^\ep_s}_{j}E^\ep_n(X)_s \right\}.
	\]
\end{example}
\begin{example}\label{ex:filt3}
	Let $A^\ep_s$ be as the above example.
	Then a closed sub-PAM $\widetilde{A}$ of $\widetilde{E}_1(X)$ can be
	given by
	\[
		\widetilde{A}=\{(\xi, \ep, s)\in \widetilde{E}_1(X)%
		~|~
		\xi\in A^\ep_s\}.
	\]
	Then we have a filtration $F^{\widetilde{A}}_j\config{n-1}{\widetilde{E}_1(X)}$.
\end{example}

\begin{rem}
	By the reason explained in the Remark in \S \ref{sec:map},
	we see that the projection map $p\co\widetilde{E}_n(X)\to \config{n-1}{\Sigma X}$
	defined in \S \ref{sec:map} preserves the filtrations given in
	Example \ref{ex:filt1} and Example \ref{ex:filt2} in the sense that
	$F_j\widetilde{E}_n(X)=p^{-1}F_j\config{n-1}{\Sigma X}$.
\end{rem}

\begin{lem} \label{thm:continuity}
	Let $A$ be a closed sub-PAM of $M$.
	Suppose that there exists a map $u\co A\to [0,1]$ such that $u^{-1}(0)=A$.
	Suppose also that
	$(M-A)\times (M-A)\cap M_2= \emptyset$ and
	$f\co M\to N$ is a function which preserves partial sums, continuous on
	$A$ and $M-A$.
	Then the induced function $C(id,f)\co \config{n}{M}\to \config{n}{N}$
	is continuous on
	$F^A_j\config{n}{M}-F^A_{j-1}\config{n}{M}$ for any $j$.
\end{lem}

\begin{proof}
	Let $\pi\co \coprod\limits_{k\geq 0} Z^{(k)}_n(M) \to \config{n}{M}$ denote the projection.
	By the definition of quotient topology, it suffices to show that
	$(id\times f)^k$ is continuous on 
	$\pi^{-1}(F^A_j\config{n}{M}-F^A_{j-1}\config{n}{M})\subset \coprod_{k\geq 0} Z_n^{(k)}(M)$.
	To do this, we express $\pi^{-1}(F^A_j\config{n}{M}-F^A_{j-1}\config{n}{M})$
	as a disjoint union
	\[
		\coprod_{k\geq j} \left(\bigcup_{[\sigma]}%
		Z^{(k)}_n(M)\cap \sigma((M-A)^j\times A^{k-j})\right)
	\]
	where  $[\sigma]$ runs over the congruence class in 
	$\Sigma_k/\Sigma_j\times \Sigma_{k-j}$
	and $\sigma\in \Sigma_k$  acts on  $M^k$  by permutation.
	Then we see that it suffices to show that $f^k\co M^k\to N^k$ is continuous on
	$\cup_{[\sigma]}\sigma((M-A)^j\times A^{k-j})$.
	By the hypothesis, $f^k$ is continuous on each $\sigma((M-A)^j\times A^{k-j})$.
	As we shall show below, each $\sigma((M-A)^j\times A^{k-j})$ is open in
	$\cup_{[\sigma]}\sigma((M-A)^j\times A^{k-j})$ and the lemma follows.
	
	To show that $\sigma((M-A)^j\times A^{k-j})$ is open, we set
	\[
		Z=\{(a_1,\dots,a_k)~|~\mbox{Max}\{u(a_i)~|~i>j\}< \mbox{Min}\{u(a_i)~|~i\leq j\}\}\subset A^k.
	\]
	Then $Z$ is an open neighborhood of $(M-A)^{j}\times A^{k-j}$ such that
	$\sigma Z\cap \sigma' Z=\emptyset$ if $[\sigma]\neq [\sigma']$ in
	$\Sigma_{k}/\Sigma_{j}\times\Sigma_{k-j}$.
	%
	%
\end{proof}

Now we prove Proposition \ref{thm:quasifibration}. The proof reduces to two
lemmas below.(We use May's form of the Dold-Thom criterion for a quasifibration
\cite{d-t},\cite{may}.)
Our proof is similar to the argument given in \S 4 of \cite{caruso-waner}, 
but we present the proof here since the construction of maps and homotopies 
are special to our setting and not obvious.
Recall that a subset $V$ of $\config{n-1}{\Sigma X}$ is said to be
distinguished if $p\co p^{-1}V\to V$ is a quasifibration.
\begin{lem}\label{thm:distinguished}
	Any open set
	$V\subset F_j\config{n-1}{\Sigma X}-F_{j-1}\config{n-1}{\Sigma X}$
	is distinguished.
\end{lem}

\begin{proof}
	We show that $p^{-1}(V)\simeq V\times \widetilde{I}_n(X)$ for any
	open set 
	$V\subset F_jC(\RR^{n-1},$ $\Sigma X)-F_{j-1}\config{n-1}{\Sigma X}$.
	Firstly, we construct a map $\varphi\co p^{-1}V\to \widetilde{I}_n(X)$.
	Suppose $\xi$ is an element of $E^\ep_1(X)_s$ represented by
	$m((J_1,x_1), \dots, (J_k,x_k))\in I^\ep_{(2k)}(X)_{(-s,s)}$.
	Let $T_t$ denote the translation of
	intervals by $t$. Then we define a function
	$\varphi^\ep_s\co E^\ep_1(X)_s\to I^\ep_1(X)_{s+2\ep}$ by
	\[
		\varphi^\ep_s(\xi)=
		\left\{
		\begin{array}{ll}
		
			[ T_{2\ep}((J_1,x_1),\dots,(J_k, x_k))] & %
			\mbox{~if~} l(J_1) \geq \ep/2 \\
			
		[(K,x_1), T_{2\ep}((J_1, x_1),\dots, (J_k, x_k)] & %
		\mbox{~if~} l(J_1) < \ep/2
		\end{array}\right. 
	\]
	where 
	$K=(\ep-l(J_1),2\ep-l(J_1),-1,-p_L(J_1))$.
	Note that the definition of $\varphi^\ep_s$ does not depend on the choice of
	the representative of $\xi$. Since $\varphi^\ep_s$ is continuous on
	$A^\ep_s$ and $E^\ep_1(X)_s-A^\ep_s$, the induced function	
	$\varphi\co E^\ep_n(X)_s\to I^\ep_n(X)_{s+2\ep}$ is continuous on
	$F^{A^\ep_s}_jE^\ep_n(X)_s-F^{A^\ep_s}_{j-1}E^\ep_n(X)_s$ by Lemma 
	\ref{thm:continuity}. Moreover $\varphi$
	induces a function
	$\widetilde{\varphi}\co \widetilde{E}_1(X)\to \widetilde{I}_1(X)$, which
	is continuous on $\widetilde{A}$ and $\widetilde{E}_1(X)-\widetilde{A}$.
	By Lemma \ref{thm:continuity}, a function
	\[
	C(id_{\RR^{n-1}},\widetilde{\varphi})\co%
	\config{n-1}{\widetilde{E}_1(X)}\to \config{n-1}{\widetilde{I}_1(X)}
	\] is continuous on
	$F^{\widetilde{A}}_j\config{n-1}{\widetilde{E}_1(X)}-%
	F^{\widetilde{A}}_{j-1}\config{n-1}{\widetilde{E}_1(X)}$.
	Note that $\widetilde{E}_n(X)$ can be embedded into
	$\config{n-1}{\widetilde{E}_1(X)}$ by the correspondence
	\[
		(\{(v_i,\xi_i)\}_i, \ep,s)\mapsto [\{(v_i,(\xi_i,\ep,s))\}_i ],
	\]
	where  $v_i\in \RR^{n-1}$ and $\xi_i\in E^\ep_1(X)_s$. Similarly,
	we can embed $\widetilde{I}_n(X)$ into $\config{n-1}{\widetilde{I}_1(X)}$.
	Then we can restrict $C(id_{\RR^{n-1}},\widetilde{\varphi})$ to a function
	$\varphi\co\widetilde{E}_n(X)\to \widetilde{I}_n(X)$,
	which is  continuous on $F_j\widetilde{E}_n(X)-F_{j-1}\widetilde{E}_n(X)$.

	Secondly we construct a map $\psi\co V\times\widetilde{I}_n(X)\to p^{-1}V$.
	Suppose $y=[t]\wedge x\in S^1\wedge X$ where  $x\in X$, $t\in [-1,1]$ and 
	$[t]\in S^1=[-1,1]/{\pm 1}$. We define a map 
	$s^\ep\co \Sigma X\to I^\ep_1(X)_{(-2\ep, 2\ep)}$ by
	\[
		s^\ep(y)=[m(L,x)] =[(-L,x),(L,x)] ,
	\]
	where  
	$L=\left(|t|\ep/2, \left(|t|/2+1\right)\ep, p ,1\right)$.
	Here, if $t\neq 0$, we put $p=-t/|t|$
	while if $t=0$ we may put either $p=+1$ or $-1$ since $s(y)$ can be
	represented by one labelled interval lying over the origin $0\in\RR$.
	Now suppose
	$v= [(c_1,y_1),\dots, (c_j,y_j)]\in V$ and $y_i\in \Sigma X$. We define
	a map $\sigma^\ep\co V\to p^{-1}V$ by
	\[
		\sigma^\ep(v)=[(c_1,s^\ep(y_1)),\dots, (c_j, s^\ep(y_j))].
	\]
	Then $\psi\co V\times \widetilde{I}_n(X)\to p^{-1}V$ is defined by
	\[
		\psi(v, (\xi,\ep,s))=%
		(\sigma^\ep(v)+m(T_{2\ep}(\xi)),\ep,s+2\ep).
	\]

	Next we show that $\psi\circ (p\times \varphi)\co p^{-1}V\to p^{-1}V$ is
	homotopic to $id_{p^{-1}V}$. Observe that $\psi\circ (p\times \varphi)$
	is induced by the function
	$\Phi^\ep_s\co E^\ep_1(X)_s\to E^\ep_1(X)_{s+4\ep}$ defined by
	\[
		\Phi^\ep_s(\xi)=
		\left\{
		\begin{array}{ll}
			[~ m( T_{4\ep}((J_1,x_1),\dots, (J_k, x_k)) )~]%
			 & \mbox{~if~} l(J_1) \geq \ep/2 \\
			
			[~ m( (L, x_1), T_{2\ep}(K, x_1), %
			 T_{4\ep}((J_1,x_1),\dots, (J_k, x_k)) ) ~]%
			& \mbox{~if~} l(J_1) < \ep/2
		\end{array}\right.\hspace{-3pt} 
	\]
	where  $K=(\ep-l(J_1), 2\ep-l(J_1), -1, -p_L(J_1))$,
	$L=(l(J_1), l(J_1)+\ep, p_L(J_1), 1)$.
	A homotopy $\psi\circ (p\times \varphi)\simeq id_{p^{-1}V}$ is induced by a deformation
	of $\Psi^\ep_s$ to $id_{E^\ep_1(X)_s}$ in $\widetilde{E}_1(X)$, where
	${E^\ep_1(X)_s}$ is regarded as a subspace of $\widetilde{E}_1(X)$ by
	the correspondence $\xi\mapsto (\xi,\ep,s)$. This deformation is given
	by a function
	\[
		H\co E^\ep_1(X)_s\times I\to \widetilde{E}_1(X)
	\]
	which coincides with $\Psi^\ep_s$ on $E^\ep_1(X)_s\times 0$ and
	is the identity on $E^\ep_1(X)_s\times 1$.
	Intuitively, the essential task of this homotopy 
	(on the right hand side of the origin) is that
	\begin{enumerate}
		\item Push $T_{2\ep}(K)$ and $T_{4\ep}((J_1),\dots, (J_k))$
		to the left until $T_*(K)$ meets $L$. Then 
		$T_*(K, x_1)$ and $(L,x_1)$ are merged.
		\item Push $T_*(J_1,\dots, J_k)$
		to the left until $T_{*'}(J_1)$ meets $L\cup T_*(K)$.
		Then $T_{*'}(J_1,x_1)$ and $(L\cup T_*(K), x_1)$ are merged.
		\item Push the right end of $L\cup T_*(K)\cup T_{*'}(J_1)$ and
		$T_{*'}((J_2,x_2),\dots, (J_k, x_k))$ to the left until
		the length of $L\cup T_*(K)\cup T_{*'}(J_1)$ coincides with
		that of original $J_1$.
	\end{enumerate}
	More precisely, we consider a homotopy $h^\ep_t\co [0,\infty)\to [0,\infty)~(0\leq t\leq 1)$
	given by the following formulae.
	
If $0\leq t\leq \frac{1}{4}$:
	\[
		h^\ep_t(u)=\left\{\begin{array}{ll}
			u & u\leq (2-4t)\ep\\
			(2-4t)\ep &  (2-4t)\ep< u\leq (2+4t)\ep\\
			u-8t\ep &  u> (2+4t)\ep
		\end{array}\right.
	\]
	If $\frac{1}{4}< t\leq \frac{1}{2}$:
	\[
		h^\ep_t(u)=\left\{\begin{array}{ll}
			u &  u\leq \ep \\
			\ep &  \ep< u\leq 3\ep\\
			u-2\ep &  3\ep < u \leq \left(\frac{9}{2}-2t\right)\ep \\
			\left( \frac{5}{2}-2t \right)\ep &	\left(\frac{9}{2}-2t\right)\ep < u 
			\leq \left( \frac{7}{2}+2t\right) \ep \\
			u-( 4t+1)\ep & \left(\frac{7}{2}+2t\right)\ep < u 
		\end{array}\right.
	\]
	If $\frac{1}{2} < t \leq \frac{3}{4}$:
	\[
		h^\ep_t(u)=\left\{\begin{array}{ll}
			u &  u\leq \ep \\
			\ep &  \ep< u\leq 3\ep\\
			u-2\ep &  3\ep < u \leq \left(\frac{9}{2}-2t\right)\ep \\
			\left(\frac{5}{2}-2t\right)\ep &
			\left(\frac{9}{2}-2t\right)\ep 
			< u \leq \left(\frac{11}{2}-2t\right) \ep \\
			u-3\ep &  u > \left(\frac{11}{2}-2t\right)\ep 
		\end{array}\right.
	\]
	If $\frac{3}{4} < t \leq 1$:
	\[
		h^\ep_t(u)=\left\{\begin{array}{ll}
			u & u \leq \left(\frac{5}{2}-2t\right) \ep\\
			\left(\frac{5}{2}-2t\right) \ep & \left(\frac{5}{2}-2t\right)\ep<u
			\leq \left(\frac{5}{2}+2t\right)\ep \\
			u-4t\ep &  u > \left(\frac{5}{2}+2t\right)\ep
		\end{array}\right.
	\]
	Then we define $H\co \widetilde{E}_1(X)\times I\to \widetilde{E}_1(X)$
	by $H( \xi, \ep ,s; t )=(h_{t*}\circ \Phi_s^\ep(\xi),\ep ,s )$.
	Since $H$ is continuous on $\widetilde{A}\times I$ and
	$\widetilde{E}_1(X)\times I-\widetilde{A}\times I,$
	we can apply Lemma \ref{thm:continuity} and get a map
	\[
		F_j^{\widetilde{A}\times I}\config{n-1}{\widetilde{E}_1(X)\times I}-%
		F_{j-1}^{\widetilde{A}\times I}\config{n-1}{\widetilde{E}_1(X)\times I}%
		\to \config{n-1}{\widetilde{E}_1(X)}.
	\]
	Consider the following sequence of embeddings
	\[
		\widetilde{E}_n(X)\times I \hookrightarrow \config{n-1}{\widetilde{E}_1(X)}\times I
		\hookrightarrow \config{n-1}{\widetilde{E}_1(X)\times I}.
	\]
	Since these embeddings are compatible with the filtrations $F_j, F^{\widetilde {A}}_j,$
	and $F^{\widetilde {A}\times I}_j,$ we can restrict $H$ to
	$(F_j\widetilde{E}_n(X)-F_{j-1}\widetilde{E}_n(X)) \times I$. 
	Since this restriction map has its image in
	$\widetilde{E}_n(X)\subset \config{n-1}{\widetilde{E}_1(X)\times I},$
	we obtain a map
	$H\co(F_j\widetilde{E}_n(X)-F_{j-1}\widetilde{E}_n(X)) \times I\to \widetilde{E}_n(X)$.
	By the definition of $h^\ep_t,$ $H$ is a fibre-preserving map with respect to
	$p\co\widetilde{E}_n(X)\to \config{n-1}{\Sigma X}$. Thus we obtain a homotopy
	$H\co p^{-1}V\times I\to p^{-1}V$ between $\psi\circ (p\times \varphi)$ and $id_{p^{-1}V}$.

	Finally we show that
	$(p\times \varphi)\circ \psi\simeq id_{V\times\widetilde{I}_n(X)}$.
	It suffices to show that $\pi_i\circ (p\times \varphi)\circ \psi \simeq \pi_i~(i=1,2),$
	where $\pi_1$ and $\pi_2$ are the projections onto $V$ and $\tilde{I}_n(X)$ respectively.
	For any $y=[t]\wedge x\in \Sigma X, t\in [-1,1], x\in X$,
	we put $t^\ep (y)=[(K,x),(L,x)]\in I^\ep_1(X)_{4\ep}$,
	where
	\[
		K=\left(\left(1-|t|/2\right)\ep, \left(2-|t|/2\right)\ep, -1, t/|t|\right), 
	\]
	$$
		L=\left(\left(2+|t|/2\right)\ep, \left(3+|t|/2\right)\ep,-t/|t|, 1\right).
	\leqno{\hbox{and}}$$
	Suppose $v=[(c_1,y_1),\dots, (c_j,y_j)]\in V, y_i=[t_i]\wedge x_i\in \Sigma X, 
	t_i\in [-1,1],$ and $x_i\in X$. We put 
	$\tau^\ep(v)=[(c_1,t^\ep(y_1)),\dots, (c_j,t^\ep(y_j))]\in I^\ep_n(X)_{4\ep}$.
	Then we have 
	$((p\times \varphi)\circ\psi)(v,(\xi,\ep,s))=(v,(\tau^\ep(v)+T_{4\ep}(\xi),\ep,s+4\ep))$.
	From this formula, it follows that
	$\pi_1\circ (p\times \varphi )\circ \psi$ coincides with $\pi_1$. 
	We prove $\pi_2\circ (p\times \varphi)\circ \psi \simeq \pi_2$ by
	the ``push to the left argument."
	Consider a homotopy $k_t\co (0,\infty]\to (0,\infty]$ defined by
	\[
		k^\ep_t(u)=\left\{\begin{array}{ll}
			u& 0\leq u\leq2\ep(1-t)\\
			2\ep(1-t) & 2\ep(1-t)< u \leq 2\ep (1+t)\\
			u-4\ep t & 2\ep(1+t)\leq u.
		\end{array}\right.
	\]
	Then we define $K\co V\times \widetilde{I}_n(X)\times I\to \widetilde{I}_n(X)$
	by
	$K( v , \xi ,\ep , s ; t ) = (k^\ep_{t*}(\xi),\ep, s)$.
	Thus we constructed a homotopy $K\co V\times \widetilde{I}_n(X)\times I\to \widetilde{I}_n(X)$
	between $\pi_2\circ (p\times \varphi)\circ \psi$ to $\pi_2$.
	This completes the proof of the lemma.
\end{proof}

\begin{lem} \label{thm:deformation}
	There exist an open neighborhood $U$ of $F_{j-1}\config{n-1}{\Sigma X}$
	 in $F_j\config{n-1}{\Sigma X}$ and homotopies $h_t\co U\to U$ and
	$H_t\co p^{-1}U\to p^{-1}U$ such that
	\begin{enumerate}
		\item $h_0=id_U$ and $h_1(U)\subset F_{j-1}\config{n-1}{\Sigma X}$,
		\item $H_0=id_{p^{-1}U}$ and $pH_t=h_tp$ for all $t,$ and
		\item $H_1\co p^{-1}z\to p^{-1}h_1(z)$ is a homotopy equivalence
		for all $z\in U$.
	\end{enumerate}
\end{lem}

\begin{proof}
	Let $u\co X\to I$ and a homotopy $k_t\co X\to X~(0\leq t\leq 1)$
	represent $(X,*)$ as a NDR-pair. Thus $u^{-1}(0)=*$ , 
	$k|_{X\times 0}=id_X$, $k(*, t)=*$ and 
	$k(W,1)=*,$ where $W=u^{-1}[0,1)$.
	We take $U\subset F_j\config{n-1}{\Sigma X}$ to be a neighborhood of
	$F_{j-1}\config{n-1}{\Sigma X}$ which consists of elements 
	represented by
	$((c_1,y_1),\dots, (c_j,y_j))$ such that there exist one or more $i$ with
	$|t_i|>\frac{1}{2}$ or $x_i\in W$, where $y_i=t_i\wedge x_i$.
	We define a homotopy $h'_t\co [-1,1]\to [-1,1], ~(0\leq t\leq 1)$ by
	\[
		h'_t(u)=\left\{\begin{array}{ll}
		-1 & -1\leq u\leq \frac{t}{2}-1\\
		\frac{2u}{2-t} & \frac{t}{2}-1\leq u\leq 1-\frac{t}{2}\\
		1 & 1-\frac{t}{2}\leq u\leq 1.
		\end{array}\right.
	\]
	Then we can define a homotopy $h_t\co U\to U$ by
	\[
		[(c_1,t_1\wedge x_1),\dots, (c_j, t_j\wedge x_j)]\mapsto%
		[(c_1,h_t(t_1)\wedge k_t(x_1)),\dots, (c_j, h_t(t_j)\wedge k_t(x_j))].
	\]
	Note that $h_0=id_U$ and $h_1(U)\subset F_{j-1}C(\RR^{n-1},\Sigma X)$.

	Next we construct a homotopy $H_t\co p^{-1}U\to p^{-1}U$ which covers $h_t$.
	Let $H'_t\co (0, 1]\to (0, 1]$ denote the homotopy defined by
	$
		H'_t(u)=\left(1-\frac{t}{2}\right)u
	$
	and $K_t=E_n(k_t)\co E_n(X)\to E_n(X)$.
	Then we can define a homotopy
	$H_t\co p^{-1}U\to p^{-1}U$ by
	$(\xi, \ep, s)\mapsto (K_t(\xi), H'_t(\ep), s)$.
	It is straightforward to check that $H_0=id_{p^{-1}U}$ and $pH_t=h_tp$ for all $t$.
	
	To show that $H_1\co p^{-1}z\to p^{-1}h_1(z)$ is a homotopy equivalence
	for all $z\in U$, a homotopy inverse map $G\co p^{-1}h_1(z)\to p^{-1}z$ is defined
	as follows:
	Suppose $z=[(c_1,y_1),\dots,(c_j,y_j)]\in U$ and $y_i=t_i\wedge x_i\in \Sigma X$.
	We put
	$g^\ep(y)=[m(K,x)]\in E^\ep_1(X)_{2\ep},$
	where $K=(|t|\ep/2, (2-|t|)\ep, -t/|t|, t/|t|),$ then
	we put $\gamma^\ep(z)=[(c_1,g^\ep(y_1)),\dots,(c_k,g^\ep(y_k))]\in E_n(X)$.
	Now $G\co p^{-1}h_1(z)\to p^{-1}z$ is defined by
	$G([\xi, \ep, s])=(\gamma^\ep(z)+T_{2\ep}(\xi), \ep, s+2\ep)$, where
	$T_{2\ep}(\xi)$ is understood to be an element of $E_n(X)$ given by
	translation $2\ep$ on the right and $-2\ep$ on the left; 
	this construction is ambiguous and not continuous as it is,
	but $G$ is well-defined and continuous by virtue of the insertion of 
	$\gamma^\ep(\xi)$.
	A proof that $G$ is the homotopy inverse of
	$H_1\co p^{-1}z\to p^{-1}h_1(z)$ is again by
	the ``push to the left" argument.
\end{proof}

%
%
%
%

%
%
%
\Addresses\recd
%
%
%
%

\begin{thebibliography}
\bibitem{caruso-waner}
\textbf{J Caruso}, \textbf{S Waner}, \emph{An approximation to {$\Omega
  \sp{n}\Sigma \sp{n}X$}}, Trans. Amer. Math. Soc. 265 (1981) 147--162
  \MR{607113}

\bibitem{cohen}
\textbf{F Cohen}, \emph{Homology of {$\Omega \sp{(n+1)}\Sigma \sp{(n+1)}X$} and
  {$C\sb{(n+1)}X,\,n>0$}}, Bull. Amer. Math. Soc. 79 (1973) 1236--1241
  (1974) \MR{0339176}

\bibitem{d-t}
\textbf{A Dold}, \textbf{R Thom}, \emph{Quasifaserungen und unendliche
  symmetrische {P}rodukte}, Ann. of Math. (2) 67 (1958) 239--281 \MR{0097062}

\bibitem{kallel}
\textbf{S Kallel}, \emph{Spaces of particles on manifolds and generalized
  {P}oincar\'e dualities}, Q. J. Math. 52 (2001) 45--70 \MR{1820902}

\bibitem{may}
\textbf{J\,P May}, \emph{The geometry of iterated loop spaces},
  Springer-Verlag, Berlin (1972) \MR{0420610}

\bibitem{mcduff}
\textbf{D McDuff}, \emph{Configuration spaces of positive and negative
  particles}, Topology 14 (1975) 91--107 \MR{0358766}

\bibitem{salvatore}
\textbf{P Salvatore}, \emph{Configuration spaces with summable labels}, from:
  ``Cohomological methods in homotopy theory (Bellaterra, 1998)'', Progr. Math.
  196, Birkh\"auser, Basel (2001)  375--395 \MR{1851264}

\bibitem{segal}
\textbf{G Segal}, \emph{Configuration-spaces and iterated loop-spaces}, Invent.
  Math. 21 (1973) 213--221 \MR{0331377}

\bibitem{shimakawa}
\textbf{K Shimakawa}, \emph{Configuration spaces with partially summable labels
  and homology theories}, Math. J. Okayama Univ. 43 (2001) 43--72 \MR{1913872}

\end{thebibliography}
\end{document}